\documentclass[11pt]{amsart}

\usepackage[margin=1in]{geometry}
\usepackage{amsmath,amssymb,amsthm,mathtools}
\usepackage{mathptmx}
\usepackage[T1]{fontenc}
\usepackage{enumitem}
\usepackage[colorlinks=true,linkcolor=blue,citecolor=blue,urlcolor=blue]{hyperref}

\numberwithin{equation}{section}
\newtheorem{theorem}{Theorem}[section]
\newtheorem{proposition}[theorem]{Proposition}
\newtheorem{lemma}[theorem]{Lemma}

\newtheorem{conjecture}[theorem]{Conjecture}
\theoremstyle{definition}

\theoremstyle{remark}

\newcommand{\N}{\mathbb{N}}
\newcommand{\R}{\mathbb{R}}
\newcommand{\Z}{\mathbb{Z}}
\newcommand{\E}{\mathbb{E}}
\renewcommand{\P}{\mathbb{P}}
\newcommand{\cC}{\mathcal{C}}
\newcommand{\cD}{\mathcal{D}}
\newcommand{\cF}{\mathcal{F}}
\newcommand{\cY}{\mathcal{Y}}
\newcommand{\one}{\mathbf{1}}
\newcommand{\hGamma}{\widehat{\Gamma}}

\title[Positive recurrence of weakly reversible networks]{Weakly reversible stochastic mass-action systems\\ are positive recurrent}

\author{Sibo Peng}
\address{Huazhong University of Science and Technology, Wuhan 430074, Hubei, China}
\email{kuhaizhongdebailanren@gmail.com}

\author{Chen Deng}
\address{Huazhong University of Science and Technology, Wuhan 430074, Hubei, China}
\email{sunnyorange.andrew@gmail.com}

\author{Yuntao Zou}
\address{Huazhong University of Science and Technology, Wuhan 430074, Hubei, China}
\email{zouyuntao@hust.edu.cn}

\author{Ruixuan Li}
\address{Huazhong University of Science and Technology, Wuhan 430074, Hubei, China}
\email{rxli@hust.edu.cn}

\date{\today}


\begin{document}

\begin{abstract}
We prove that every finite weakly reversible reaction network, modelled as a continuous-time Markov chain with stochastic mass-action kinetics, is positive recurrent on each of its communicating classes for every choice of positive rate constants. The minimal chain is non-explosive, every state has a finite mean return time, and on each infinite class the stationary distribution $\pi$ has factorial-type tails: for every $\delta>0$ there is $C_\delta<\infty$ with $x!\,\pi(x)\le C_\delta\exp\big(\delta\, m(x)\log(2+m(x))\big)$, where $m(x)=\max_i x_i$. This resolves a conjecture of Anderson and Kim. The proof constructs a Lyapunov function for the embedded jump chain from local potentials, one for each logarithmic direction in which the class is unbounded. Each local potential is defined through minimal-peak levels of the network projected onto its direction and is controlled by a layered Poisson sampling argument. The local potentials are glued by an affine approximation lemma with arbitrary tolerances, and the resulting shell drift, combined with a pointwise return estimate, gives a factorial bound on the occupation measure of an excursion. The factorial occupation bound and the finiteness of the mean return time have been formally verified in Lean~4.
\end{abstract}

\maketitle

\tableofcontents

\section{Introduction}\label{sec:intro}

Reaction networks are used throughout biology and chemistry to describe interacting populations of molecules. When the molecules are abundant, the concentrations are modelled by a system of ordinary differential equations. When some species are present in low numbers, as is common inside cells, the natural model is a continuous-time Markov chain on $\N^d$ whose state records the copy numbers of the $d$ species \cite{Kur72,AK15}. A recurring theme of reaction network theory is to deduce dynamical properties from the reaction graph alone, for every choice of rate constants. The deficiency zero theorem of Horn, Jackson and Feinberg \cite{HJ72} is the prototype in the deterministic setting. Its stochastic counterpart, due to Anderson, Craciun and Kurtz \cite{ACK10}, shows that weakly reversible networks of deficiency zero have product-form Poisson stationary distributions, and are therefore positive recurrent, for all rate constants.

A network is \emph{weakly reversible} if every reaction $y\to y'$ lies on a directed cycle of the reaction graph. For the deterministic model, weak reversibility is the hypothesis of the global attractor and persistence conjectures, and of the existence theorem of Boros \cite{Bor19}. For the stochastic model, Anderson and Kim conjectured that weak reversibility alone suffices for positive recurrence, calling this the positive recurrence conjecture \cite[\S6]{AK18}. The question is also listed as Problem~272 in the collection \emph{AIM --- Open Applied Problems} \cite{AIM}.

\begin{conjecture}[Anderson--Kim positive recurrence conjecture]\label{conj:AK}
Every weakly reversible reaction network with stochastic mass-action kinetics is positive recurrent on each of its communicating classes, for every choice of positive rate constants.
\end{conjecture}

The conjecture is delicate because intuition from the deterministic model fails in both directions. There are networks whose deterministic model has a globally stable positive equilibrium while the stochastic model goes extinct almost surely \cite{AEJ14}. There are also strongly endotactic networks, whose deterministic model is permanent \cite{GMS14}, for which the stochastic model is transient and even explosive \cite{ACKN20}. Partial results are known. Positive recurrence holds for complex-balanced networks \cite{ACK10}; for weakly reversible networks with a single linkage class under additional structural assumptions \cite{ACK20}; for a subclass of strongly endotactic networks identified through tier sequences \cite{ACKN20}; and for weakly reversible networks with one-dimensional stoichiometric subspace \cite{WX20}. Non-explosivity is known for complex-balanced networks \cite{ACKK18} and for bimolecular weakly reversible networks \cite{Xu24}. In general, every communicating class of a weakly reversible network is closed \cite{PCK14}. Closedness alone does not exclude transience, since an infinite closed class can be transient; the open question is whether every such class is positive recurrent, that is, whether returns occur almost surely with finite mean return time.

Our main result settles Conjecture~\ref{conj:AK}.

\begin{theorem}\label{thm:main}
Let $G$ be a finite weakly reversible reaction network on $d$ species and let $\kappa_r>0$ be arbitrary rate constants. Consider the stochastic mass-action model \eqref{eq:rates} on $\N^d$. Then:
\begin{enumerate}[label=\textup{(\roman*)}]
\item every communicating class is closed;
\item on every communicating class the minimal chain is non-explosive, and every non-absorbing state has a finite mean first return time; in particular the chain is positive recurrent on every class, and absorbing states carry point-mass stationary distributions;
\item on every infinite communicating class $\Gamma$, the unique stationary distribution $\pi$ satisfies, for every $\delta>0$,
\begin{equation}\label{eq:main-tail}
x!\,\pi(x)\le C_\delta\exp\big(\delta\,m(x)\log(2+m(x))\big)\qquad(x\in\Gamma),
\end{equation}
where $x!=\prod_i x_i!$, $m(x)=\max_ix_i$, and $C_\delta<\infty$ depends on $\delta$, the network, the rate constants and $\Gamma$.
\end{enumerate}
\end{theorem}

The bound \eqref{eq:main-tail} says that $\pi(x)$ is at most $1/x!$ up to a factor $e^{o(m\log m)}$. This matches the decay of the product-form Poisson distributions of complex-balanced networks \cite{ACK10}, so in this sense the stationary tails of every weakly reversible network are no heavier than those of complex-balanced networks, even when the stationary distribution is not of product form.

\subsection*{Strategy of the proof}
The difficulty is that in a general weakly reversible network different reactions dominate in different regions of $\N^d$, and the region boundaries are curved in the original coordinates but nearly polyhedral in logarithmic coordinates $w_R(x)_i=\log(x_i+1)/\log R$. No single explicit function decreases on average everywhere. We therefore build a Lyapunov function for the embedded jump chain in three stages.

\emph{Local potentials.} For each direction $a\in[0,1]^d$ in which the communicating class is unbounded (a \emph{usable} direction, Section~\ref{sec:directions}), we project the network onto the coordinates where $a$ vanishes and attach to each state a real \emph{level} $\ell_a(x)=a\cdot x+\psi_a(x_I)$. Here $\psi_a$ is the minimal peak needed to descend arbitrarily far in the projected network (Section~\ref{sec:witness}). In a logarithmic neighbourhood of $a$, after a time change, all reaction rates become powers of $R$ whose exponents (the \emph{relative orders}) take finitely many values. A layered Poisson sampling argument (Section~\ref{sec:sampling}) shows that the level drops within a time polynomial in $R$ before a positive-order reaction can raise it. This gives a local potential $U_{a,R}$ comparable to $R^{\theta a\cdot x}$ that solves an exact Poisson equation for the jump chain (Section~\ref{sec:local}).

\emph{Gluing.} The local potentials are combined into
$V_R=\sum_\lambda R^{\theta(b_\lambda(R)+a_\lambda\cdot\one)}U_{a_\lambda,R}$
using a finite family of affine functions whose maximum $F_R$ has the following property: any face that is within $(\log R)^2$ of the maximum has its slope within a prescribed tolerance of the current logarithmic position (Section~\ref{sec:affine}). This lemma, with tolerances that may depend discontinuously on the slope, lets us switch off every local potential outside its own domain at a cost that is negligible against the gain of the active one. The result is the shell drift estimate $(P-I)V_R\le -R^{-K}R^{\theta F_R}$ on shells $uN\le m(x)<N$ (Section~\ref{sec:drift}).

\emph{Occupation bound.} A pointwise return estimate based on a coordinatewise descent lemma controls the number of entrances into a shell. Combined with the shell drift, it gives the bound $x!L_o(x)\le C_\delta e^{\delta m\log m}$ for the occupation measure of an excursion from a root $o$ (Section~\ref{sec:occupation}). Summing over $x$ gives positive recurrence and \eqref{eq:main-tail} (Section~\ref{sec:pr}).

The proof uses only the strong Markov property, Dickson's lemma, Tonelli's theorem, and elementary properties of Poisson clocks. The layered construction in the affine approximation lemma of Section~\ref{sec:affine} follows an idea of \cite[\S2]{OAI}; we give a complete proof, adapted to the margin $(\log R)^2$ required here.

\subsection*{Organisation}
Section~\ref{sec:model} sets up the model and notation. Sections~\ref{sec:directions}--\ref{sec:levels} develop the geometry of communicating classes, bounded descent witnesses and levels. Section~\ref{sec:sampling} contains the layered sampling lemma and Section~\ref{sec:local} the local potentials. Section~\ref{sec:affine} proves the affine approximation lemma. Section~\ref{sec:drift} glues the potentials and proves the shell drift. Section~\ref{sec:occupation} derives the factorial occupation bound and Section~\ref{sec:pr} completes the proof of Theorem~\ref{thm:main}. Section~\ref{sec:remarks} discusses the constants, the scope of the result and formal verification.

\section{The model}\label{sec:model}

Let $d\ge1$ be the number of species and let the reactions be indexed by a finite set. Reaction $r$ has source $s_r\in\N^d$, target $t_r\in\N^d$ and increment $\Delta_r=t_r-s_r$; we keep only reactions with $s_r\neq t_r$. If $x\ge s_r$ coordinatewise, the chain may jump from $x$ to $x+\Delta_r$ along label $r$ at rate
\begin{equation}\label{eq:rates}
q_r(x)=\kappa_r\prod_i\frac{x_i!}{(x_i-(s_r)_i)!},\qquad \kappa_r>0,
\end{equation}
and $q_r(x)=0$ otherwise. Write $q(x)=\sum_rq_r(x)$, $m(x)=\max_ix_i$, $x!=\prod_ix_i!$ and $J_*=\max_r\|\Delta_r\|_\infty$. When there is at least one reaction, put $\kappa_*=\min_r\kappa_r>0$ and $p=\max_r\sum_i(s_r)_i$. Every enabled rate is at least $\kappa_*$, and the total rate is at most $C(1+m(x))^p$.

The \emph{complexes} are the vectors $s_r,t_r$; the reaction graph has the complexes as vertices and an edge $s_r\to t_r$ for each $r$. Weak reversibility means that for every reaction $s_r\to t_r$ there is a directed word in the reaction graph from $t_r$ back to $s_r$.

\begin{lemma}\label{lem:closed}
Every communicating class of the chain is closed.
\end{lemma}

\begin{proof}
Let $r$ be enabled at $x$ and put $z=x-s_r\ge0$, the \emph{remainder}. Replace each vertex $c$ of a return word from $t_r$ to $s_r$ in the reaction graph by the actual state $z+c$. This gives a feasible path from $x+\Delta_r$ back to $x$, since each reaction on it is enabled at $z+c\ge c$. We call the set of states $z+c$ with this fixed remainder, together with the reactions between them, a \emph{complex copy}. Thus the two endpoints of every enabled edge communicate, and every communicating class is closed.
\end{proof}

Finite non-absorbing classes are handled by the return theorem for finite irreducible chains, and absorbing states carry point masses. From now on we fix an arbitrary infinite class $\Gamma$ and a root $o\in\Gamma$. Then $q(x)\ge\kappa_*$ on $\Gamma$.

Let $q(x,y)$ be the sum of the rates of labels jumping from $x$ to $y$, let $P(x,y)=q(x,y)/q(x)$ be the transition matrix of the jump chain, and let $Q(x,y)=q(x,y)$ for $y\ne x$, $Q(x,x)=-q(x)$. Let $(X_n)$ be the jump chain, $J_o=\inf\{n\ge1:X_n=o\}$ (possibly infinite), and define for the chain started at $o$
\begin{equation}\label{eq:occupation}
L_o(x)=\E_o\sum_{n<J_o}\one_{\{X_n=x\}},\qquad G(x)=\frac{L_o(x)}{q(x)}.
\end{equation}
Attaching to each jump an exponential holding time with mean $1/q(X_n)$ constructs the continuous-time chain. Even before non-explosion is known, Tonelli's theorem identifies $G(x)$ as the mean physical time spent at $x$ before the first return to $o$. We shall prove that for every $\delta>0$
\begin{equation}\label{eq:factorial}
x!\,L_o(x)\le C_\delta\exp\big(\delta\,m(x)\log(2+m(x))\big)\qquad(x\in\Gamma).
\end{equation}
Multiplying \eqref{eq:factorial} by $1/q(x)\le1/\kappa_*$ gives the same bound for $x!G(x)$. Since $x!\ge m(x)!$ and the number of lattice points of radius $r$ is polynomial in $r$, any fixed $0<\delta<1$ makes $L_o$ summable. Section~\ref{sec:pr} shows how this yields positive recurrence in real time.

\section{Directions of a communicating class}\label{sec:directions}

For $a\in[0,1]^d$ let $I(a)=\{i:a_i=0\}$. We call $a$ \emph{usable} if there is $v\in\N^d\setminus\{0\}$ with $v_{I(a)}=0$ and $o+v\in\Gamma$. In words, some nonzero growth vector of the class is supported in the positive coordinates of $a$.

\begin{lemma}\label{lem:unusable}
If $a$ is not usable, there is a finite constant $C_a$ such that every $x\in\Gamma$ satisfies
\begin{equation}\label{eq:unusable}
m(x)\le C_a\Big(1+\max_{i\in I(a)}x_i\Big),
\end{equation}
where the maximum over the empty set is $0$.
\end{lemma}

\begin{proof}
Let $M=\{v\in\N^d:o+v\in\Gamma\}$. A feasible word of labels can be translated by any nonnegative vector and remains feasible, because the enabling conditions are monotone. If $u,v\in M$, translating a path $o\to o+v$ by $u$ shows $u+v\in M$. If moreover $u\ge v$, translating a path $o+v\to o$ by $u-v$ leads from $o+u$ to $o+u-v$, so $u-v\in M$. Hence $D=M-M$ is an additive subgroup of $\Z^d$ and
\begin{equation}\label{eq:M}
M=D\cap\N^d.
\end{equation}
By Dickson's lemma, $\Gamma$ has finitely many coordinatewise minimal elements; call this set $B_0$. For $x\in\Gamma$ choose $b\in B_0$ with $b\le x$; translating a path $b\to o$ by $x-b$ gives $x-b\in M$. Conversely, if $v\in M$, translating a path $o\to b$ by $v$ gives $b+v\in\Gamma$. Therefore $\Gamma=\bigcup_{b\in B_0}(b+M)$.

Suppose \eqref{eq:unusable} fails. Then there is a sequence $x^{(n)}\in\Gamma$ with $m_n=m(x^{(n)})\to\infty$ and $\max_{i\in I(a)}x^{(n)}_i/m_n\to0$. Pass to a subsequence with $x^{(n)}/m_n\to v_*$ and with a fixed minimal element $b$ below $x^{(n)}$. Since $\|v_*\|_\infty=1$, the limit is nonzero; its coordinates are nonnegative and vanish on $I(a)$. As $(x^{(n)}-b)/m_n$ lies in the real span of $D$, so does $v_*$.

Let $U=\{i:(v_*)_i>0\}$. The intersection of the real span of $D$ with the coordinate subspace $\{v:v_{U^c}=0\}$ is defined by rational linear equations. Approximating $v_*$ inside this intersection by rational vectors keeps strict positivity on $U$. Choosing a real basis consisting of integer vectors of $D$, the rational approximant has rational coefficients in this basis. Multiplying by a suitable positive integer yields a nonzero $v\in D\cap\N^d$ supported in $U$. By \eqref{eq:M}, $o+v\in\Gamma$, contradicting the assumption that $a$ is not usable.
\end{proof}

From now on we use the scale $R=2N$ and the logarithmic coordinates
\[
w_R(x)_i=\frac{\log(x_i+1)}{\log R}.
\]
If $a$ is not usable and $\|w_R(x)-a\|_\infty<e$ for some $e<1$, then $x_i+1<R^e$ on $I(a)$, and \eqref{eq:unusable} gives $m(x)\le C_a(1+R^e)$. Hence, for any fixed $u>0$ and all large $N$, such states lie outside the shell $m(x)\ge uN$.

\section{Projected words and bounded descent witnesses}\label{sec:witness}

\subsection{Phases and copy heights}
Fix a direction $a$ and write $I=I(a)$. A \emph{phase} $y\in\N^I$ records the inventory on the coordinates $I$; in addition we keep a real \emph{height} $H$. In this auxiliary graph, a label $r$ may be used when $y\ge(s_r)_I$; it changes $y$ to $y+(\Delta_r)_I$ and $H$ to $H+a\cdot\Delta_r$. The \emph{copy height} of this edge is
\begin{equation}\label{eq:copyheight}
h=H-a\cdot s_r,\qquad S=\max_{c\in\cC}a\cdot c,
\end{equation}
where $\cC$ is the set of complexes. After the edge, the new height is $h+a\cdot t_r$. Fix the remainder phase $y-(s_r)_I$ and the scalar $h$, and place each complex $c$ at phase $y-(s_r)_I+c_I$ and height $h+a\cdot c$. This gives the same \emph{projected complex copy}. All edges in this copy have the same copy height, and the return words of the reaction graph are return words in the auxiliary graph. Feasibility is preserved by increasing the phase inventory, by lowering the initial height, and by raising the permitted upper bound on copy heights.

The height increment of a word is the sum of $a\cdot\Delta_r$ over its labels. A word is \emph{feasible under the cap $b$} if every step has the required phase inventory and every edge has copy height at most $b$. A nonempty feasible word started at height zero requires $b\ge-S$.

\subsection{The bounded witness lemma}

\begin{lemma}\label{lem:witness}
Fix $C\ge0$ and $D_0>0$. There is an integer $B(a,C,D_0,I)$ with the following property. For every phase $y$, every $b\in[-S,C]$ and every $0<u\le D_0$, if some word started at $(y,0)$ is feasible under the cap $b$ and reaches height $\le-u$, then such a word exists with length at most $B(a,C,D_0,I)$.
\end{lemma}

\begin{proof}
We induct on the number of phase coordinates. During the induction we allow an arbitrary subset of phase coordinates, always keeping the same scalar increments and the copy height formula \eqref{eq:copyheight}.

If there are no phase coordinates, a descending word must contain an edge with negative increment. Before the first such edge the height is nonnegative, so this edge can already be used from height zero. Afterwards we repeat this edge; each repetition starts lower, so remains feasible. The minimum absolute value of the finitely many negative increments is positive, so the number of steps needed to reach $-D_0$ is uniformly bounded. If there is no edge with negative increment, the hypothesis is void.

Suppose the set of phase coordinates is nonempty. Let $C'=C+D_0$ and $D'=D_0+C+S$. Forget one phase coordinate at a time and take the inductive bound in the lower dimension with parameters $C',D'$; let $B_1$ be the maximum of these bounds. Let $s_*=\max_{r,i}(s_r)_i$ and $K=1+B_1s_*$. If some phase coordinate has inventory at least $K$, we may first take a word of length at most $B_1$ in the system that forgets this coordinate. When the coordinate is restored, the consumption before the $j$-th edge is at most $(j-1)s_*$, so the inventory suffices to enable every edge and the height trajectory is unchanged. This handles a large initial inventory and also the exit points described below.

Now assume all initial phase coordinates are below $K$. Scan the given word in order, keeping a modified prefix until the target is first reached, the phase first leaves the cube, or a repeatable descending loop appears. At every scanning time there is a suffix that is still feasible and reaches the target. The phases in the kept prefix lie in the finite cube $\{0,\dots,K-1\}^I$ and are distinct. If the cube has $M_0$ points, the prefix has at most $M_0$ edges. Consider the next edge.

If the new height is already $\le-u$, the prefix is a bounded witness. If the new phase is in the cube and appears for the first time, append it to the prefix. If the new phase equals a previously kept phase, the end of the prefix forms a closed phase segment of length at most $M_0$; let $\Delta H$ be its height increment. If $\Delta H\ge0$, delete the segment and continue with the remaining labels from the first occurrence of this phase. The phase conditions are the same and every height of the suffix decreases by $\Delta H$, so the cap and the terminal condition still hold. The number of unprocessed labels decreases, so deletion can be repeated.

If $\Delta H<0$, repeat the closed segment from its starting point. The phase returns to itself and each round starts lower, so each round respects the original cap. Since the segment has length at most $M_0$, its total increment takes values in a finite set; the negative values have absolute value bounded below by some $\rho>0$. The starting height of the segment is at most $C+S$, since it is either zero or the landing height of an edge that is already feasible. Hence at most $1+\lceil(C+S+D_0)/\rho\rceil$ rounds reach the target. Together with the kept prefix before the segment, the length is still uniformly bounded.

Finally suppose the new phase leaves the cube for the first time without reaching the target, and let $H'$ be the height at this point. Then $-u<H'\le b+S$. Resetting the height to zero at this point, the remaining descent target is $u'=u+H'>0$ and the relative cap is $b'=b-H'$. These satisfy
\begin{equation}\label{eq:witness-params}
0<u'\le D_0+C+S,\qquad -S\le b'\le C+D_0.
\end{equation}
At the exit point some inventory is at least $K$. The projection of the remaining feasible word satisfies the hypotheses of the induction in one dimension lower, so the suffix can be replaced by a word of length at most $B_1$, which is lifted back to the original dimension by the inventory argument above. The original word is finite and each scanning step processes at least one new label. The cases above exhaust the ways in which the scan can end, and the maximum length over all branches gives the required constant.
\end{proof}

\subsection{Minimal peaks and levels}
Fix a usable direction $a$. Let $\cY_a$ be the set of phases from which there is a closed phase word with strictly negative increment. Such a word can be translated together with the initial phase inventory, so $\cY_a$ is upward closed. Every projected edge can be undone by a path inside the same copy, whose total height increment is the negative of that of the edge. So one can return to the original phase, run the negative loop, and go forward again; this attaches a negative loop to any adjacent phase. Hence $\cY_a$ is invariant under projected edges.

By usability, choose a nonzero $v$ with $v_I=0$ and $o+v\in\Gamma$. A word from $o+v$ to $o$ projects to a word with the same initial and final phase and increment $-a\cdot v<0$. A word from any state of $\Gamma$ to $o$ can also be projected and undone edge by edge. Hence the phase of every state of $\Gamma$ lies in $\cY_a$. Define
\[
\hGamma_a=\{x\in\N^d:x_I\in\cY_a\}.
\]
It contains $\Gamma$ and is closed under reactions and under nonnegative translations; it is used below to construct translated witnesses.

By Dickson's lemma, $\cY_a$ has finitely many minimal phases. For each, choose a negative loop started at height zero, and let $C_a\ge0$ bound the maximal copy height of all these loops. Translating the word of a minimal phase upward gives the same bound at every phase $y\in\cY_a$. Repeating the negative loop only lowers the copy heights, so it reaches arbitrarily low heights.

Fix $D_a>C_a+S+1$. Any witness descending from height zero to height $\le-D_a$ ends at a phase that is still in $\cY_a$. Appending a negative loop there, the added copy heights are at most $-D_a+C_a<-S$. Therefore the smallest cap under which one can descend arbitrarily low equals the smallest cap under which one first descends by $D_a$, and this minimum lies in $[-S,C_a]$. Indeed, for any feasible cap, Lemma~\ref{lem:witness} lets us replace the descending word by a word of bounded length. Started at height zero, the copy height of every edge is determined by the word of labels. Finitely many such words give finitely many candidate peaks, and every phase has at least one feasible candidate.

Denote this minimal peak by $\psi_a(y)$. It takes finitely many values. If the phase inventory increases, the old witness can be translated, so $\psi_a$ is nonincreasing. The \emph{level} of an actual state is
\begin{equation}\label{eq:level}
\ell_a(x)=a\cdot x+\psi_a(x_I).
\end{equation}
Equivalently, $\ell_a(x)\le b$ if and only if from $(x_I,a\cdot x)$ one can reach arbitrarily low heights under the cap $b$. The minimal-peak witnesses from a state $x$ have uniformly bounded length. Each of their edges has copy height at most $k=\ell_a(x)$, and they land at height at most $a\cdot x-D_a$, so the landing level is at most $a\cdot x-D_a+C_a<k$. If the initial inventory of every coordinate with $a_i>0$ exceeds the bound on the witness length times the maximal source inventory, each step of the witness is also enabled in that coordinate, and we obtain an actual descending word.

\section{Level preservation on complete copies and relative orders}\label{sec:levels}

Fix a usable $a$ and drop the subscript when no confusion arises. Let $x\to x'$ be an enabled edge with copy height $h$. Returning from $x'$ to $x$ along the same copy and then following any arbitrarily-low descending word from $x$ shows $\ell(x')\le\max(\ell(x),h)$. Following the edge forward from $x$ and then a descending word from $x'$, and reversing the roles, gives $\ell(x)\le\max(\ell(x'),h)$.

If $k=\ell(x)$ and $h<k$, then $\ell(x')=k$, and by the same inequalities every point of the return path stays at level $k$. If $h=k$, every point of the copy has level at most $k$, and the return path may first enter a lower level. If a jump raises the level, then necessarily $h>k$, and the new level is at most $h$.

Within the level $\ell=k$, define the \emph{relative order} of a label by
\begin{equation}\label{eq:relorder}
d_r(x,k)=h-k=-a\cdot s_r-\psi_a(x_I).
\end{equation}
It belongs to a finite set $\cD_a$ independent of $x$ and $k$. The change of level along any actual jump also lies in a fixed finite set, because it equals $a\cdot\Delta_r+\psi_a(x'_I)-\psi_a(x_I)$. Hence any decrease of level is at least some fixed $g_a>0$. If the level increases along a label of relative order $d'$, the increase is at most $d'$.

A path that uses only labels strictly faster than a given nonpositive order has a return path inside the original level that uses only faster edges. For each step of a general path of nonpositive order, the return path either goes back along the same relative order or first reaches a lower level. These are the two forms of reversibility used in the sampling induction below.

\section{Layered Poisson sampling}\label{sec:sampling}

\subsection{Time change and Poisson clocks}\label{sec:clocks}
Fix a usable direction $a$ and $\eta>0$. Let $\Omega_{a,R}$ be the set of states of $\hGamma_a$ with $m(x)\le R-2$ and $\|w_R(x)-a\|_\infty<E_a$. Choose
\[
E_a<\tfrac12\min_{a_i>0}a_i .
\]
Then for $R$ large, inside this domain all inventories in the positive coordinates of $a$ are large enough to lift the descending witnesses of Section~\ref{sec:witness}. We also require $pE_a\le\eta$, where $p$ is the maximal source order from Section~\ref{sec:model}.

For $n\ge s\ge1$ and $0\le j<s$ we have $(n-j)/(n+1)\ge1/(s+1)$. This gives a two-sided comparison, up to fixed factors, between $q_r(x)$ and $\prod_i(x_i+1)^{(s_r)_i}$ for enabled labels, allowing zero inventories and zero source orders. Multiplying all rates out of a state $x$ by the same factor $R^{\ell(x)-a\cdot x}$ does not change the jump chain. On the level $\ell(x)=k$ the resulting rates satisfy
\begin{equation}\label{eq:timechanged}
c_0R^{-d'-\eta}\le\widehat q_r(x)\le C_0R^{-d'+\eta},\qquad d'=d_r(x,k).
\end{equation}
Here we used $|(w_R(x)-a)\cdot s_r|\le pE_a$ and $k-a\cdot x+a\cdot s_r=-d'$. All constants can be chosen uniformly over the finitely many labels.

For each relative order $d'$ we run an independent Poisson clock of rate $\Lambda_{d'}=K_0R^{-d'+\eta}$. When this clock rings, each available label of this order is selected with probability $\widehat q_r(x)/\Lambda_{d'}$, and with the remaining probability the ring is void. Taking $K_0\ge2C_0$ times the number of labels ensures that the void probability is at least $1/2$ and that each specified enabled label is selected with probability at least $c_1R^{-2\eta}$. Labels that are not enabled receive no probability. As long as the process stays in the domain and in the same level, this thinning realises exactly the true rates. When the region or the level changes, the thinning continues with the level of the current state.

Discard the positive orders and number the nonpositive orders $d_1<\dots<d_L=0$; a clock of order zero that always rings void may be added. Write $\Lambda_j=\Lambda_{d_j}$. Since the order differences in a fixed direction are finite, for large $R$ we have $2\Lambda_j\le\Lambda_{j-1}$ for every $j>1$. The thinned process using the first $j$ orders is stopped when it reaches $\ell<k$ or leaves $\Omega_{a,R}$; the stopped state may be recorded as an absorbing mark.

\subsection{Static targets and a uniform word length}
Let $A(x,k)$ denote a target predicate. We require that it contains all states with $\ell(x)<k$ and that it is preserved by nonnegative translations $(x,k)\mapsto(x+v,k+a\cdot v)$. Static reachability is defined in the whole of $\hGamma_a$. A path may use, on level $k$, edges whose order index is at most $j$, and it ends when it meets a lower level. If such a path leaves the domain first in the actual thinned process, the stopping event has already occurred, so static witnesses need not be restricted to the domain.

Fix $k-a\cdot x=\sigma$. For starting points on level $k$, $\sigma$ can only be one of the finitely many values of $\psi_a$. For each $\sigma$, consider the starting points from which the target can be reached using the first $j$ orders. Dickson's lemma gives finitely many minimal points, and every such starting point dominates one of them coordinatewise. Choose a target witness for each minimal point. When the witness is translated to a given starting point, the enabling conditions are preserved. Since $\psi_a$ is nonincreasing, the level along the way is at most the translated level $k+a\cdot v$. If the level drops early, the truncated word reaches the target. If it stays at the original level, the copy heights and the level both increase by $a\cdot v$, so the relative orders are unchanged. The translated original endpoint also lies in the target. The maximum length over the witnesses of all minimal points gives a uniform integer $B=B(a,j,A)$, fixed before $\eta$ and the scale $R$ are chosen.

Let $m_j(x,k)$ be the minimum number of order-$j$ edges over all words that reach the target using the first $j$ orders, with $m_j=\infty$ if there is none. Reachable states have $m_j\le B$. If a segment of a fast word uses only edges with index less than $j$, its return path also stays on the original level and uses only fast edges. Concatenating the return path with a target word shows that fast words preserve the finite budget $m_j$.

General allowed edges also preserve \emph{target reachability}. For a landing point that is still alive, follow the return path inside the same copy. If it first reaches a lower level, that level belongs to the target. If it stays on level $k$, returning to the starting point and following the original target word suffices. If in the actual process the state has already left the domain or dropped a level, there is nothing to prove. This argument allows arbitrarily many allowed edges before sampling.

For $b\ge1$ define the \emph{entrance target} $A_b$ to consist of all lower-level states, together with the points still on level $k$ that have an enabled order-$j$ edge leading to a state with $m_j\le b-1$. An edge that already descends also qualifies. $A_b$ is preserved by nonnegative translations: translate the edge and the subsequent budget witness; if the level drops early we have succeeded, and otherwise the relative orders and the budget are preserved. If $m_j(x,k)=b>0$, take a witness achieving the minimum count. Before its first order-$j$ edge it uses only fast edges, so $A_b$ can be reached by fast edges. After crossing this edge the budget is at most $b-1$.

\subsection{The sampling estimate given the number of rings}

\begin{lemma}[Sampling lemma]\label{lem:sampling}
Consider a target as above and a starting point on the level from which the target is reachable. Run the thinned process using the first $j$ orders and sample it at an independent $\mathrm{Exp}(\lambda)$ time with $0<\lambda\le\Lambda_j$. Then the probability that at the sampling time the process has stopped or is in the target is at least $c_{j,A}R^{-b_{j,A}\eta}$. The positive constant $c_{j,A}$ and the finite exponent $b_{j,A}$ are uniform in the starting point, the level and the admissible sampling rate; the threshold for $R$ may depend on $\eta$.
\end{lemma}

\begin{proof}
We induct on the number of available orders. If no order is available, reachability of the target means that the starting point is already in the target, so the claim holds for every sampling rate. Assume the claim for the faster orders. Let the order-$j$ clock compete with the sampling clock. The time of the next ring and its type are independent: the time is $\mathrm{Exp}(\Lambda_j+\lambda)$ and the type is order $j$ with probability $p_j=\Lambda_j/(\Lambda_j+\lambda)\ge1/2$. For instance, the joint density that the order-$j$ clock rings first at a time in $dt$ is $\Lambda_j\exp(-(\Lambda_j+\lambda)t)\,dt$, and factorising it gives the independence. After an order-$j$ ring we continue by the memoryless property of exponential clocks until the sampling clock rings.

Let $n$ be the number of order-$j$ rings before sampling. The probability of at least $B$ rings is $p_j^B\ge2^{-B}$. Conditionally on exactly $n\ge B$ rings, the order of the two kinds of rings is fixed. The $n+1$ gaps are still independent $\mathrm{Exp}(\Lambda_j+\lambda)$ times, independent of the faster clocks and of the random marks used for thinning. In each gap the faster process can be restarted from its state at the beginning of the gap and the induction hypothesis applied. If $j>1$, the induction allows this rate because $\Lambda_j+\lambda\le2\Lambda_j\le\Lambda_{j-1}$. If $j=1$, we use the base case with no available orders.

Let the first $n-B$ rings happen arbitrarily. As long as the process has not stopped, target reachability is preserved, so $m_j\le B$ before the last $B$ rings begin. Suppose $b$ controlled rings remain. If $m_j\le b-1$, the path in this fast gap preserves the budget, and we choose the ring to be void. If $m_j=b$, we use the fast sampling induction to reach $A_b$ or to have stopped by the end of the gap. If the process is still alive, we select the corresponding enabled label, which lowers the budget to $b-1$. In both branches the conditional success probability of each round is at least $c_*R^{-b_*\eta}$, where we take a common lower bound over the finitely many bounds for $A_1,\dots,A_B$, the void ring and single labels; $b_*$ can be taken as the maximum of the faster exponents plus $2$.

After these $B$ rounds, if the process has not stopped, the budget is zero, which means the target can be reached using only faster orders. The last gap is the fast gap before the sampling clock, and we apply the fast sampling induction once more to the target itself. Multiplying the conditional probabilities shows that, given $n\ge B$, the success probability is at least
\begin{equation}\label{eq:sampling-product}
\big(c_*R^{-b_*\eta}\big)^B\,c_{j-1,A}R^{-b_{j-1,A}\eta}.
\end{equation}
This lower bound does not depend on the extra $n-B$ rings. Removing the conditioning with $\P(n\ge B)\ge2^{-B}$ gives the claim. Because of the uniform word length, the induction uses only a fixed number of entrance targets. If $B=0$, the product over the $B$ rings is empty and only the last fast gap is used. The constants and exponents are chosen first from the graph and the target; afterwards $\eta$ and the logarithmic tolerance are tightened.
\end{proof}

\section{Local potentials}\label{sec:local}

\subsection{Mean stopping time of the nonpositive-order process}
On a given level $k$ take as target the set of all lower-level states. By Section~\ref{sec:witness}, inside the domain there is an actual feasible descending word, each of whose edges has copy height at most $k$; if a prefix of it has already left the domain, stopping has also occurred. Let $S_k$ be the first time at which the process using only the clocks of nonpositive order drops a level or leaves the domain. Apply the sampling lemma at an independent observation clock of rate $\Lambda_L$. For every state that is still alive, the conditional probability of stopping before the next observation is at least $p_R=cR^{-b\eta}$.

Let $0=T_0<T_1<\cdots$ be the independent observation times. By the strong Markov property and the independent increments of the observation clock, $\P(S_k>T_i)\le(1-p_R)^i$. The mean $1/\Lambda_L$ of a future observation gap is independent of the event of survival up to $T_i$. Summing the interval lengths over these events gives
\begin{equation}\label{eq:stoptime}
\E S_k\le\frac{1}{\Lambda_Lp_R}\le CR^{b\eta}.
\end{equation}
The last inequality deliberately keeps a generous exponent. It is uniform in the initial level and the starting point in the domain.

\subsection{Restoring positive orders and controlling the update weight}
Add the neglected positive-order clocks back, independently. Starting from a point on a level, call a \emph{segment} the period up to a level drop, an exit from the domain, or the first ring of a positive-order clock, whichever comes first. Let $\zeta$ be the duration of a segment. When a positive order rings first, the reaction is selected at that ring according to \eqref{eq:timechanged}, and the segment ends even if the ring is void. The next segment starts from the actual state left behind. In every state the actual label rates of the whole process are the time-changed mass-action rates, so its sequence of actual jumps is still governed by $P$.

Given a trajectory with only nonpositive orders, the probability that the first ring of a positive order $d'>0$ precedes $S_k$ is $1-\exp(-\Lambda_{d'}S_k)$. By $1-\exp(-t)\le t$ and \eqref{eq:stoptime}, this clock wins with probability at most
\begin{equation}\label{eq:preempt}
\Lambda_{d'}\E S_k\le CR^{-d'+(b+1)\eta}.
\end{equation}
If a segment ends with a level drop inside the domain, the level falls by at least $g_a$. If it ends at a positive order $d'$ inside the domain, the level rises by at most $d'$; this includes void rings and level drops. After an exit from the domain the weight is set to zero. For fixed $0<\theta<1$, the conditional expectation of the surviving weight at the end of a segment, divided by the weight at its start, is therefore at most
\begin{equation}\label{eq:weightratio}
R^{-\theta g_a}+\sum_{d'>0}CR^{-(1-\theta)d'+(b+1)\eta}.
\end{equation}
If positive orders exist, first take $(b+1)\eta<\tfrac12(1-\theta)\min_{d'>0}d'$; if there are none the sum is empty. Then for $R$ large, \eqref{eq:weightratio} is at most $1/2$. The exponent $b$ of the sampling lemma was fixed before choosing $\eta$. This fixes, for each $a$ and $\theta$, a tolerance $E_a>0$ for which all the conditions of Section~\ref{sec:sampling} hold.

\subsection{From clock rewards to a potential for the jump chain}
Every segment's actual jumps start on the same level $k$. The number of actual jumps from the start to the end of a segment is at most the total number of rings of all orders, whose expectation is $\E\zeta\sum_{d'}\Lambda_{d'}$. This follows by stopping the martingale ``Poisson count minus rate times time'' at $\zeta\wedge n$, and then using \eqref{eq:stoptime} and monotone convergence. Since $\zeta\le S_k$ and there are finitely many fixed orders, the expected number of jumps is at most a fixed power of $R$. Each actual jump from level $k$ pays $R^{\theta k}$, so the conditional expected reward of a segment is at most $CR^{C}R^{\theta k}$.

Iterating \eqref{eq:weightratio} from segment to segment, the expected surviving initial weight of the $i$-th segment is at most $2^{-i}R^{\theta\ell(x)}$. Sum the nonnegative rewards with Tonelli's theorem. With $\tau_\Omega=\inf\{n\ge0:X_n\notin\Omega_{a,R}\}$ define
\begin{equation}\label{eq:localpotential}
U_{a,R}(x)=\E_x\sum_{n<\tau_\Omega}R^{\theta\ell_a(X_n)},\qquad U_{a,R}(x)=0\quad(x\notin\Omega_{a,R}).
\end{equation}
The terms of the sum correspond exactly to the initial states of the actual jumps before leaving the domain. The total time-changed rate of every state in the domain is positive and finite, so the first actual jump occurs at a finite time almost surely. Comparing the first term of \eqref{eq:localpotential} with the geometric sum above, and using that $\psi_a$ is bounded, we obtain a fixed $K_a$ such that for all large $R$
\begin{equation}\label{eq:localcompare}
R^{-K_a}R^{\theta a\cdot x}\le U_{a,R}(x)\le R^{K_a}R^{\theta a\cdot x}\qquad(x\in\Omega_{a,R}).
\end{equation}
The first-step decomposition also gives the exact Poisson equation for the jump chain
\begin{equation}\label{eq:poisson}
(P-I)U_{a,R}(x)=-R^{\theta\ell_a(x)}\qquad(x\in\Omega_{a,R}).
\end{equation}
Here and below we use the extension by zero, so jumps that leave the logarithmic domain in one step are included in the identity.

\section{Affine approximation with a margin}\label{sec:affine}

The layered construction in this section follows an idea of \cite[\S2]{OAI}. We give a complete proof, in the near-active form with the margin $(\log R)^2$ that is needed in Section~\ref{sec:drift}.

\subsection{Statement}
Let $E:[0,1]^d\to(0,\infty)$ be an arbitrary positive tolerance function, not necessarily continuous. We construct finitely many slopes $a_\lambda\in[0,1]^d$ and offsets $b_\lambda(R)=o(R)$ such that, with
\begin{equation}\label{eq:affine-family}
L_{\lambda,R}(x)=a_\lambda\cdot(x+\one)+b_\lambda(R),\qquad F_R(x)=\max_\lambda L_{\lambda,R}(x),
\end{equation}
there is $R_0$ for which, whenever $R\ge R_0$ and $w_R(x)\in[0,1]^d$ (here $x$ may be any real vector with $x_i>-1$ and $w_R(x)_i=\log(x_i+1)/\log R$),
\begin{equation}\label{eq:near-active}
L_{\lambda,R}(x)\ge F_R(x)-(\log R)^2\ \Longrightarrow\ \|w_R(x)-a_\lambda\|_\infty<E(a_\lambda).
\end{equation}

\emph{Minimum form.} Put $h=1/R$, $p=-w_R(x)\in[-1,0]^d$, $r=-a$, $z=x+\one$ (so $z_i=h^{p_i}$), $c=-b$ and $g(h)=(\log(1/h))^2=(\log R)^2$. Then $L=-(r\cdot z+c)$ and $F_R=-\min(r\cdot z+c)$, and \eqref{eq:near-active} says: every face within $g(h)$ of the minimum satisfies $\|p-r\|_\infty<E_-(r):=E(-r)$. The offset condition $b=o(R)$ becomes $c=o(h^{-1})$. It therefore suffices to prove Proposition~\ref{prop:Pd} below with $A=-1$, $B=0$.

\subsection{The inductive statement}
For $A<B\le0$, $0<h<1$ and $p\in\R^d$ write $Z_h(p)=(h^{p_1},\dots,h^{p_d})$. A \emph{face family} is a finite set of pairs $(r,c_r)$ with $r\in[A,B]^d$ and $c_r:(0,1)\to\R$. The value of the face $r$ at $(h,p)$ is $V_r(h,p)=r\cdot Z_h(p)+c_r(h)$. The face $r$ is \emph{near-active} at $(h,p)$ if $V_r(h,p)\le\min_{r'}V_{r'}(h,p)+g(h)$.

\begin{proposition}\label{prop:Pd}
For every $A<B\le0$ and every positive tolerance $E_-:[A,B]^d\to(0,\infty)$ there is a face family $\cF$ with offsets $c_r(h)=o(h^A)$ as $h\to0$ and the following \emph{sequential property}:

\textup{(S)} If $h_n\to0$, $p^{(n)}\to p\in[A,B]^d$ \textup{(}the $p^{(n)}$ themselves need not lie in the cube\textup{)}, and a fixed face $r\in\cF$ is near-active at $(h_n,p^{(n)})$ for all $n$, then $\|p-r\|_\infty<E_-(r)$.
\end{proposition}

\emph{The sequential property implies the uniform property.} Assume (S). We claim there is $h_0$ such that for all $h<h_0$ and all $p\in[A,B]^d$, every near-active face satisfies $\|p-r\|_\infty<E_-(r)$. Otherwise there are $h_n\to0$, $p^{(n)}\in[A,B]^d$ and near-active faces $r_n$ that violate the tolerance. The family is finite, so a subsequence has $r_n\equiv r$; the cube is compact, so a further subsequence has $p^{(n)}\to p$. The violation $\|p^{(n)}-r\|_\infty\ge E_-(r)$ persists in the limit, contradicting (S). Taking $A=-1$, $B=0$ gives \eqref{eq:affine-family}--\eqref{eq:near-active} with $R_0=1/h_0$.

We prove Proposition~\ref{prop:Pd} by induction on $d$. Assume it holds in every dimension $k<d$ (for $d=1$ no hypothesis is needed), and fix $A<B\le0$ and $E_-$.

\subsection{Local families with a fixed limiting minimal exponent}
Fix $t\in[A,B]$. The local family $\cF_t$ consists of faces of types $0,1,\dots$, constructed as follows.

\emph{Type $0$.} One face, with slope $t\one$ and offset $0$. If $d=1$ or $t=B$, the local family consists of this face only.

\emph{Type $k$ ($1\le k\le d-1$, so that $t<B$).} Proceed for $k=1,\dots,d-1$ in turn.
\begin{enumerate}[label=(\arabic*)]
\item Choose $\beta_k$ with $t<\beta_k<B$, $\beta_k-t<\min\{E_-(r):r\text{ of type }k-1\}$, and, if $k>1$, $\beta_k<\beta_{k-1}$. The right-hand side of the second condition is a minimum of finitely many positive numbers, so such $\beta_k$ exists.
\item For each $k$-element set $U\subset\{1,\dots,d\}$, embed a slope $\rho\in[\beta_k,B]^U$ as $\iota_U\rho\in[A,B]^d$, equal to $\rho$ on $U$ and to $t$ off $U$. On $[\beta_k,B]^U$ define the tolerance $E_U(\rho)=E_-(\iota_U\rho)$.
\item By the induction hypothesis in dimension $k<d$, obtain a face family $\cF_U$ with offsets $d_\rho(h)=o(h^{\beta_k})$.
\item The faces of type $k$ are $\iota_U\rho$, for all $k$-element sets $U$ and $\rho\in\cF_U$, with offsets $c(h)=-h^{\beta_k}+d_\rho(h)$.
\end{enumerate}
Since $\beta_k>t$ we have $h^{\beta_k}=o(h^t)$, so all offsets of the local family are $o(h^t)$. Since $\beta_k<B\le0$ we have $g(h)=o(h^{\beta_k})$; this is used repeatedly below.

\begin{lemma}\label{lem:local-seq}
Let $h_n\to0$, $p^{(n)}\to p\in[A,B]^d$ with $\min_ip_i=t$, and let a fixed face $r\in\cF_t$ be near-active in $\cF_t$ at $(h_n,p^{(n)})$ for all $n$. Then $\|p-r\|_\infty<E_-(r)$.
\end{lemma}

\begin{proof}
Write $z^{(n)}=Z_{h_n}(p^{(n)})$. We estimate $|p_i-r_i|$ coordinate by coordinate.

\emph{(i) Raised positions of a face of type $k\ge1$.} Let $U$ be the raised set. Compare with the type-$0$ face: $V_r\le V_{t\one}+g$. The two faces agree off $U$, so
\begin{equation}\label{eq:typek-vs-type0}
\sum_{i\in U}(r_i-t)\,h_n^{p^{(n)}_i}\le h_n^{\beta_k}-d_\rho(h_n)+g(h_n)=(1+o(1))\,h_n^{\beta_k}.
\end{equation}
Each coefficient satisfies $r_i-t\ge\beta_k-t>0$ and each term is nonnegative. If $p_i<\beta_k$ for some $i\in U$, then $h_n^{p^{(n)}_i}/h_n^{\beta_k}\to\infty$, contradicting \eqref{eq:typek-vs-type0}. Hence $p_U\in[\beta_k,B]^U$.

Now compare faces with the same $U$. They share the terms off $U$ and the term $-h^{\beta_k}$, so their differences equal the differences of the corresponding faces of $\cF_U$. The minimum of the local family is at most the minimum over its $\cF_U$ part, so the face $\rho$ is also near-active in $\cF_U$ along $(h_n,p^{(n)}_U)$, whose limit $p_U$ lies in $[\beta_k,B]^U$. By the sequential property of $\cF_U$,
\[
\|p_U-\rho\|_\infty<E_U(\rho)=E_-(r).
\]

\emph{(ii) Positions $j\notin U$} (for type $0$, $U=\emptyset$). Since $t=\min_ip_i$, we have $p_j\ge t=r_j$; we must show $p_j-t<E_-(r)$.
\begin{itemize}
\item If $k\le d-2$ (including $k=0$, in which case $d\ge2$), we claim $p_j\le\beta_{k+1}$. Suppose $p_j>\beta_{k+1}$ and let $U'=U\cup\{j\}$. On $U$, part (i) gives $p_i\ge\beta_k>\beta_{k+1}$ (for $k=0$, $U$ is empty). So all limiting exponents on $U'$ are strictly greater than $\beta_{k+1}$. Take any $\rho'\in\cF_{U'}$ and the type-$(k+1)$ face $r'=\iota_{U'}\rho'$. Since $h_n^{p^{(n)}_i}=o(h_n^{\beta_{k+1}})$ for $i\in U'$,
\[
V_{r'}-V_{t\one}=\sum_{i\in U'}(r'_i-t)h_n^{p^{(n)}_i}-h_n^{\beta_{k+1}}+d_{\rho'}=-(1+o(1))h_n^{\beta_{k+1}}.
\]
The original face satisfies $V_r-V_{t\one}\ge-h_n^{\beta_k}+d_\rho=-(1+o(1))h_n^{\beta_k}$ if $k\ge1$, and $V_r-V_{t\one}=0$ if $k=0$. Since $\beta_{k+1}<\beta_k$, $h^{\beta_k}=o(h^{\beta_{k+1}})$, so $V_r-V_{r'}\ge(1+o(1))h_n^{\beta_{k+1}}$. As $g=o(h^{\beta_{k+1}})$, this contradicts the near-activity of $r$. Thus $0\le p_j-t\le\beta_{k+1}-t$, and by the choice of $\beta_{k+1}$ the right-hand side is strictly less than the tolerance of every face of type $k$, in particular $E_-(r)$.
\item If $k=d-1$, the unique unraised position $j$ has $p_j=t$: all other positions have $p_i\ge\beta_{d-1}>t$, while $\min_ip_i=t$.
\item If $d=1$ or $t=B$, only type $0$ exists. For $d=1$, $p=t$; for $t=B$, $p\in[A,B]^d$ has minimum $B$, so $p=B\one$.
\end{itemize}
Combining (i) and (ii) gives $\|p-r\|_\infty<E_-(r)$.
\end{proof}

\subsection{Extending a local family to a neighbourhood}

\begin{lemma}\label{lem:local-nbhd}
For each $t\in[A,B]$ there is $\rho_t>0$ such that the conclusion of Lemma~\ref{lem:local-seq} holds for all sequences whose limit satisfies $|\min_ip_i-t|<\rho_t$. Moreover $\rho_t$ can be chosen with $t+\rho_t<\min_k\beta_k$ \textup{(}this condition is void if the local family has only type $0$\textup{)}.
\end{lemma}

\begin{proof}
For each face $r\in\cF_t$ let $\Lambda_r$ be the set of limits $p\in[A,B]^d$ for which there exist $h_n\to0$ and $p^{(n)}\to p$ such that $r$ is near-active in $\cF_t$ at every $(h_n,p^{(n)})$. The set $\Lambda_r$ is closed: if $p^{[m]}\in\Lambda_r$ and $p^{[m]}\to p$, choose from a sequence realising $p^{[m]}$ one term $(h,p')$ with $h<1/m$ and $|p'-p^{[m]}|<1/m$; the diagonal sequence realises $p$. Let
\[
K_t=\bigcup_{r\in\cF_t}\Big(\Lambda_r\cap\{p\in[A,B]^d:\|p-r\|_\infty\ge E_-(r)\}\Big).
\]
It is a finite union of compact sets, hence compact, and by Lemma~\ref{lem:local-seq} it contains no point with $\min_ip_i=t$. Since $\min_ip_i$ is continuous and does not take the value $t$ on $K_t$, there is $\rho_t>0$ with $|\min_ip_i-t|\ge\rho_t$ on $K_t$ (any $\rho_t$ works if $K_t$ is empty). Shrinking $\rho_t$ does not affect this, so we may also require $t+\rho_t<\min_k\beta_k$.
\end{proof}

Write $I_t=(t-\rho_t,t+\rho_t)$. Because $t+\rho_t<\min_k\beta_k$, all offsets of the local family are $o(h^{t+\rho_t})$: type $0$ has offset zero, and the offsets of type $k$ are of order $h^{\beta_k}$.

\subsection{A chain of intervals}
The local families are glued along a chain of the intervals $I_t$. We construct the chain by a supremum argument.

\begin{lemma}\label{lem:chain}
Let $\rho:[A,B]\to(0,\infty)$ be any function and $I_t=(t-\rho_t,t+\rho_t)$. There are $m\ge1$, centres $t_1,\dots,t_m\in[A,B]$ and cut points $A=q_0<q_1<\dots<q_m=B$ such that
\begin{enumerate}[label=\textup{(\roman*)}]
\item $q_{j-1}\in I_{t_j}$ and $q_j\in I_{t_j}$ for $1\le j\le m$, so that $[q_{j-1},q_j]\subset I_{t_j}$;
\item $t_j<q_j<t_{j+1}$ for $1\le j<m$.
\end{enumerate}
\end{lemma}

\begin{proof}
For $s\in(A,B]$, a \emph{chain to $s$} is a choice of $m$, $t_1,\dots,t_m$ and $q_0,\dots,q_m$ satisfying all the requirements of the lemma except that $q_m=s$ instead of $q_m=B$. In every chain $q_{m-1}\le t_m$: for $m=1$ this is $A\le t_1$, and for $m\ge2$ it is part of (ii). Two operations produce new chains.
\begin{enumerate}[label=(\Alph*)]
\item \emph{Moving the last cut point.} If $s'>q_{m-1}$ and $s'\in I_{t_m}$, replacing $q_m$ by $s'$ turns a chain to $s$ into a chain to $s'$.
\item \emph{Appending a centre.} Let $s'\in[A,B]$ and let $q'$ satisfy $\max(t_m,\,s'-\rho_{s'})<q'<\min(t_m+\rho_{t_m},\,s')$. Replacing $q_m$ by $q'$ and adding $t_{m+1}=s'$, $q_{m+1}=s'$ turns a chain to $s$ into a chain to $s'$. Indeed $q'\in I_{t_m}$, $q'>t_m\ge q_{m-1}$, $q'\in I_{s'}$, $t_m<q'<t_{m+1}$, and $s'\in I_{s'}$.
\end{enumerate}
Let $S\subset(A,B]$ be the set of points $s$ to which a chain exists. It is nonempty: $s_0=\min(B,A+\rho_A/2)$ belongs to $S$, with $m=1$ and $t_1=A$. Let $s^*=\sup S\in(A,B]$, and choose $s\in S$ with $s>s^*-\rho_{s^*}$ together with a chain to $s$. Since $s=q_m\in I_{t_m}$, we have $t_m-\rho_{t_m}<s<t_m+\rho_{t_m}$.

If $s^*<t_m+\rho_{t_m}$, then $s^*\in I_{t_m}$ and $s^*\ge s>q_{m-1}$, so (A) gives a chain to $s^*$. Otherwise $s^*\ge t_m+\rho_{t_m}$. The interval $\big(\max(t_m,s^*-\rho_{s^*}),\,t_m+\rho_{t_m}\big)$ is then nonempty, because $s^*-\rho_{s^*}<s<t_m+\rho_{t_m}$, and its points are less than $s^*$; so (B) gives a chain to $s^*$. In both cases $s^*\in S$.

Finally, suppose $s^*<B$ and take a chain to $s^*$, with last centre $t_m$. Then $s^*<t_m+\rho_{t_m}$, and (A) with $s'=\min\big(B,(s^*+t_m+\rho_{t_m})/2\big)>s^*$ gives a chain to a point of $S$ larger than $\sup S$, which is impossible. Hence $s^*=B\in S$.
\end{proof}

Lemma~\ref{lem:chain} replaces the usual argument of taking a finite subcover and sorting it, and it is the form of the argument checked in Lean. Note that the cut points $q_1,\dots,q_{m-1}$ are strictly less than $B\le0$, since $q_j<t_{j+1}\le B$, and that $q_\ell>A$ for $\ell\ge1$.

\subsection{Gluing into a global family}
Apply Lemma~\ref{lem:chain} with the radii $\rho_t$ of Lemma~\ref{lem:local-nbhd}. The global family $\cF$ consists of the faces of all $\cF_{t_j}$, $1\le j\le m$, where every offset of the $j$-th family is decreased by
\[
\sigma_j(h)=\sum_{\ell=1}^{j-1}h^{q_\ell}.
\]
The sum starts at $\ell=1$. Since $q_\ell>A$, $\sigma_j=o(h^A)$, and the local offsets are $o(h^{t_j})\subset o(h^A)$, so all global offsets are $o(h^A)$.

\begin{proof}[Proof of the sequential property \textup{(S)} for $\cF$]
Let a fixed face $r$ from the $j$-th family be near-active in $\cF$ along $(h_n,p^{(n)})$ with limit $p$, and let $\tau=\min_ip_i$.

\emph{(a) If $j<m$ and $\tau>q_j$, we get a contradiction.} Compare with the type-$0$ face $r'=t_{j+1}\one$ of the $(j+1)$-st family, whose offset is $-\sigma_{j+1}=-\sigma_j-h^{q_j}$. Then
\[
V_r-V_{r'}=(r-r')\cdot z^{(n)}+c_r(h_n)+h_n^{q_j}.
\]
All $p_i\ge\tau>q_j$, so $z^{(n)}_i=o(h_n^{q_j})$. Also $c_r=o(h^{t_j+\rho_{t_j}})$ and $q_j<t_j+\rho_{t_j}$ (because $q_j\in I_{t_j}$), so $c_r=o(h^{q_j})$. Hence $V_r-V_{r'}=(1+o(1))h_n^{q_j}$. Since $q_j<0$, this is much larger than $g(h_n)$, contradicting near-activity.

\emph{(b) If $j>1$ and $\tau<q_{j-1}$, we get a contradiction.} Compare with the type-$0$ face $r'=t_{j-1}\one$ of the $(j-1)$-st family, whose offset is $-\sigma_{j-1}=-\sigma_j+h^{q_{j-1}}$. Then
\[
V_r-V_{r'}=\sum_i(r_i-t_{j-1})z^{(n)}_i+c_r(h_n)-h_n^{q_{j-1}}.
\]
All slopes of the $j$-th family have coordinates $\ge t_j$, so each coefficient satisfies $r_i-t_{j-1}\ge t_j-t_{j-1}>0$ and each term is nonnegative. For a coordinate $i$ with $p_i=\tau$ we have $z^{(n)}_i/h_n^{q_{j-1}}\to\infty$. Moreover $c_r=o(h^{q_{j-1}})$, because $q_{j-1}\in I_{t_j}$ gives $q_{j-1}<t_j+\rho_{t_j}$. Hence $V_r-V_{r'}\ge(t_j-t_{j-1})z^{(n)}_i(1-o(1))$, which is much larger than $g(h_n)$, a contradiction.

\emph{(c) Conclusion.} By (a) and (b), $\tau\in[q_{j-1},q_j]\subset I_{t_j}$; for $j=1$ the lower bound $A$ and for $j=m$ the upper bound $B$ hold automatically. All faces of the $j$-th family are shifted by the same $-\sigma_j$, so near-activity in $\cF$ implies near-activity in $\cF_{t_j}$, since the global minimum is at most the minimum over the $j$-th family. By Lemma~\ref{lem:local-nbhd}, $\|p-r\|_\infty<E_-(r)$.
\end{proof}

This proves Proposition~\ref{prop:Pd} and hence \eqref{eq:affine-family}--\eqref{eq:near-active}. The slopes depend only on the tolerance function $E$ and not on $R$; the offsets depend only on $R$. Throughout the paper the slope family is fixed first and $R$ is then taken large.

\section{Glued potential and shell drift}\label{sec:drift}

Fix $\theta\in(0,1)$ as above and any $0<e<1/4$. For each usable direction let $E(a)$ be the smaller of the tolerance of Section~\ref{sec:local} and $e$; for directions that are not usable take $E(a)\le e$. Apply Section~\ref{sec:affine} to obtain a finite face family and fix its slopes and offsets. Fix any $0<u<1$, take an integer $N$ and $R=2N$, and let
\begin{equation}\label{eq:shell}
A_{u,N}=\{x\in\Gamma:\lceil uN\rceil\le m(x)<N\}.
\end{equation}
For all states with $\lceil uN\rceil-J_*\le m(x)\le N+J_*$ and all large $N$ we have $m(x)\le R-2$ and $w_R(x)\in[0,1]^d$. In this range every near-active face is usable. Indeed, if a near-active face had an unusable direction, \eqref{eq:near-active} and \eqref{eq:unusable} would give $m(x)\le C_a(1+R^e)$, contradicting $m(x)\ge uN-J_*$ for large $N$; since the face family is finite this exclusion is uniform. Every near-active usable face also satisfies the tolerance that defines its own logarithmic domain.

Let $\Lambda_*$ be the index set of the usable faces and define
\begin{equation}\label{eq:glued}
V_R(x)=\sum_{\lambda\in\Lambda_*}R^{\theta(b_\lambda(R)+a_\lambda\cdot\one)}\,U_{a_\lambda,R}(x),
\end{equation}
where each local potential is extended by zero as in \eqref{eq:localpotential}. By \eqref{eq:localcompare}, every nonzero weighted term equals $R^{\theta L_{\lambda,R}(x)}$ up to a factor $R^{\pm K}$. An active term exists and its domain contains $x$. Hence on the shells above and their one-step landing points there is a uniform constant $K$ with
\begin{equation}\label{eq:Vcompare}
R^{-K}R^{\theta F_R(x)}\le V_R(x)\le R^{K}R^{\theta F_R(x)}.
\end{equation}

We now compute the drift at $x\in A_{u,N}$. For the terms whose logarithmic domain contains $x$, \eqref{eq:poisson} gives the negative weighted contribution $-R^{\theta\ell_a(x)}$. Taking one active term, the absolute value of the total negative drift of these terms is at least $R^{-K_1}R^{\theta F_R(x)}$; here we use that $\ell_a(x)-a\cdot x=\psi_a(x_I)$ is uniformly bounded.

For a term whose domain does not contain $x$, we have $U_{a,R}(x)=0$. A positive contribution can only come from actual landing points $x'$ reached in one step with $U_{a,R}(x')>0$. Since $x$ lies within a cutoff of uniform size and the logarithmic condition of this term fails at $x$, \eqref{eq:near-active} gives $L_{\lambda,R}(x)<F_R(x)-(\log R)^2$. One jump changes any affine value by at most $dJ_*$. Estimating the term at the landing point by \eqref{eq:localcompare}, and using that the transition probabilities sum to at most one, the total positive contribution of the terms outside their domains is at most
\begin{equation}\label{eq:outside}
R^{K_2+\theta dJ_*-\theta(\log R)^2}R^{\theta F_R(x)}.
\end{equation}
For fixed $K_1,K_2$ and large $R$, \eqref{eq:outside} is less than half of the lower bound for the active negative term. Increasing $K$, we obtain uniformly
\begin{equation}\label{eq:shelldrift}
(P-I)V_R(x)\le-R^{-K}R^{\theta F_R(x)}\qquad(x\in A_{u,N}).
\end{equation}

Two endpoint estimates are also needed. Let $b_*(R)=\max_\lambda|b_\lambda(R)|=o(R)$. All slopes lie in $[0,1]^d$, so every point satisfies $F_R(x)\le d(m(x)+1)+b_*(R)$. If an outer landing point satisfies $N\le m(x)\le N+J_*$, take a maximal coordinate $i$ and an active slope $a$. Since $w_R(x)_i\ge\log(N+1)/\log R=1-o(1)$, we have $a_i\ge1-e-o(1)$. The other slope coordinates are nonnegative, so
\begin{equation}\label{eq:outerF}
F_R(x)\ge(1-e-o(1))N-b_*(R)\ge N/2.
\end{equation}
Together with \eqref{eq:Vcompare} this bounds the potential below at actual outer landing points. The constant $K$ is determined by the finite local families; the shell ratio $u$ only affects how large $N$ must be.

\section{A factorial bound on the occupation of an excursion}\label{sec:occupation}

\subsection{A pointwise visit bound without assuming return}

\begin{lemma}[Coordinatewise descent]\label{lem:descent}
There is a finite constant $C$ such that for every $x\in\Gamma$ there is a feasible path from $x$ to $o$ of length at most $C(1+m(x))$, every state $y$ of which satisfies $y_i\le x_i+C$ for all $i$.
\end{lemma}

\begin{proof}
On the whole lattice let $\mathcal U$ be the set of states $x$ from which some finite feasible word leads to a state $x'\le x$ with $\sum_ix'_i<\sum_ix_i$. Enabling of mass-action reactions is monotone under coordinatewise translation, so $\mathcal U$ is upward closed. By Dickson's lemma it has finitely many minimal elements. For each minimal element $u$ fix such a descending word. If the current state satisfies $x\in\mathcal U$, choose a minimal $u\le x$ and translate its word by $x-u$. The new state is coordinatewise no larger than $x$ and its total is smaller by at least one. The lengths of the fixed words and the amounts by which their intermediate states exceed their starting points have a common bound. After at most $\sum_ix_i\le d\,m(x)$ repetitions the state lies outside $\mathcal U$. Since each segment starts from the endpoint of the previous one, which is coordinatewise at most $x$, the excess over the original coordinates does not accumulate.

In the closed irreducible class $\Gamma$, a state outside $\mathcal U$ is a coordinatewise minimal element of $\Gamma$: if $\Gamma$ contained a strictly smaller comparable state, irreducibility would give a path to it. There are finitely many such minimal elements. Fix a path from each of them to $o$ and absorb the lengths and coordinates of these paths into the constant.
\end{proof}

By Lemma~\ref{lem:descent}, from every $x\ne o$ there is a feasible path to $o$ of length at most $C(1+m(x))$ whose intermediate radii are at most $m(x)+C$. Deleting the part before the last visit to $x$, the path does not revisit its starting point before reaching the root. On this path, each enabled label is the next actual jump with probability at least $\kappa_*/[C(1+m(x)+C)^p]$. Multiplying these probabilities, there is $C_1$ independent of $x$ such that from $x$ the chain reaches the root before returning to $x$ with probability at least $1/M_x$, where
\begin{equation}\label{eq:Mx}
M_x=\exp\big(C_1(1+m(x))\log(2+m(x))\big).
\end{equation}
From any starting point, the mean number of visits to a non-root point $x$ before the first hit of the root is at most $M_x$. Indeed, after the first hit of $x$, each visit has probability at least $1/M_x$ of being followed by a hit of the root before the next visit. By the strong Markov property the number of visits has a geometric tail with parameter $1/M_x$; paths that never revisit $x$ only decrease the count. Hence, also within a positive-return excursion from the root, $L_o(x)\le M_x$ and $L_o(o)=1$. At this point the excursion is still allowed to have infinitely many jumps.

\subsection{Stopping at actual entrances and exits of a shell}
Fix $u$ and $N$ so that the root lies inside. Let $l_N=\lceil uN\rceil$ and $\tau_N=\inf\{n\ge0:m(X_n)\ge N\}$, and consider the jump chain up to $J_o\wedge\tau_N$. For large $N$, the gap between the inner and outer boundaries exceeds the maximal jump size. List the successive parts of the excursion that lie in $l_N\le m<N$. Each part starts at an actual entrance point from the inside and stops at the actual landing point of the first return to $m<l_N$ or of the first exit to $m\ge N$. Only parts that start before $J_o\wedge\tau_N$ are counted.

The state from which an entrance is made satisfies $l_N-J_*\le m<l_N$, and each visit produces at most one entrance. The expected number of entrances is therefore at most the sum of $M_y$ over this thin shell, with $1$ for the root if it is included. The thin shell has at most $O(R^d)$ points, and for these points
\begin{equation}\label{eq:logMy}
\log M_y\le C_1(uN+O(1))\log R.
\end{equation}
An actual entrance point $z$ satisfies $m(z)\le l_N+J_*$. By \eqref{eq:Vcompare} and $F_R(z)\le d(m(z)+1)+o(N)$, every entrance point has $V_R(z)\le R^{duN+o(N)}$. So the expected sum of the initial potentials of all parts is at most
\begin{equation}\label{eq:entrancebudget}
R^{C_2uN+o(N)},\qquad C_2=C_1+d.
\end{equation}
The $o(N)$ is taken after fixing all parameters and absorbs the number of shell points and polynomial factors. The leading coefficient $C_2$ does not depend on $\theta$, $e$ or $u$.

Within a part started at the entrance point $z$, let $\rho$ be the jump index at which it leaves the shell, and let $c_R(x)=R^{-K}R^{\theta F_R(x)}$. Summing the one-step drift \eqref{eq:shelldrift}, stopped at $\rho\wedge n$, gives
\begin{equation}\label{eq:stoppeddrift}
\E_zV_R(X_{\rho\wedge n})+\E_z\sum_{j<\rho\wedge n}c_R(X_j)\le V_R(z).
\end{equation}
For a finite horizon only a finite region and its one-step neighbours are involved, so all terms are integrable. In the first term keep only the terminal potential of parts that have exited through the outside; the second term is the accumulated reward. Both increase with the horizon, so the horizon can be sent to infinity. Conditioning on the actual starting point of each part and summing over the parts with Tonelli's theorem, the right-hand side is controlled by \eqref{eq:entrancebudget}.

If $\tau_N<J_o$, exactly one part ends at an outer landing point, where $V_R\ge R^{-K+\theta N/2}$ by \eqref{eq:Vcompare} and \eqref{eq:outerF}. Hence the escape probability $p_N=\P_o(\tau_N<J_o)$ satisfies
\begin{equation}\label{eq:escape}
p_N\le R^{C_2uN+o(N)-\theta N/2}.
\end{equation}
For any $x\in A_{u,N}$, every visit to $x$ before escape or return to the root pays at least $c_R(x)$. Absorbing fixed powers of $R$, the mean number $L^{<\tau_N}_o(x)$ of such truncated visits satisfies
\begin{equation}\label{eq:truncatedvisits}
L^{<\tau_N}_o(x)\le R^{C_2uN+o(N)-\theta F_R(x)}.
\end{equation}

\subsection{Choosing the scales and combining early and late visits}
Fix $\delta>0$. First choose $\theta$ close to $1$ and $e$ close to $0$ with $d(1-\theta)<\delta/8$, $de<\delta/8$ and $e<1/4$. Choose $A\ge2$ with
\begin{equation}\label{eq:chooseA}
\theta A/4>C_1+d+2.
\end{equation}
Then take $u>0$ so small that $uA<1/2$, $C_2u<\theta/8$ and $C_2uA<\delta/8$. For a point with large $r=m(x)$ take $N=\lceil Ar\rceil$ and $R=2N$. Then $x\in A_{u,N}$. Since all parameters are fixed, the $o(N)$ in \eqref{eq:escape} is at most $\theta N/8$ at large scales, so $p_N\le R^{-\theta N/4}$.

The mean number of visits to $x$ after the first escape is at most $p_NM_x$; this follows from the strong Markov property at the actual landing point on $\{\tau_N<J_o\}$ and the pointwise bound of Section~\ref{sec:occupation} for arbitrary starting points. Using $\log x!\le dr\log(2+r)$ and \eqref{eq:chooseA},
\begin{equation}\label{eq:late}
x!\,p_NM_x\le\exp\big((d+C_1-\theta A/4+o(1))\,r\log(2+r)\big)\le1.
\end{equation}
The visits before escape are controlled through an active slope. Take a face $a$ active at $x$. The tolerance bound and $\log x_i!\le(x_i+1)\log(x_i+1)$ give
\begin{equation}\label{eq:Flower}
F_R(x)\log R\ge\log x!-ed(r+1)\log R-b_*(R)\log R.
\end{equation}
Insert \eqref{eq:Flower} into \eqref{eq:truncatedvisits}. Since $N=Ar+O(1)$, $\log R=\log(2+r)+O(1)$ and $b_*(R)=o(N)=o(r)$, we obtain
\begin{equation}\label{eq:early}
x!\,L^{<\tau_N}_o(x)\le\exp\big((d(1-\theta)+\theta ed+C_2uA+o(1))\,r\log(2+r)\big).
\end{equation}
The three fixed small quantities sum to less than $3\delta/8$. Combining \eqref{eq:late} and \eqref{eq:early} gives \eqref{eq:factorial} for large $r$. The finitely many remaining states are absorbed into $C_\delta$ using \eqref{eq:Mx}. We also record the statement for the physical occupation:
\begin{equation}\label{eq:physical}
\forall\delta>0\ \exists C_\delta<\infty\ \forall x\in\Gamma:\quad x!\,G(x)\le C_\delta\exp\big(\delta\,m(x)\log(2+m(x))\big).
\end{equation}

\section{Positive recurrence and the stationary distribution}\label{sec:pr}

\begin{proof}[Proof of Theorem~\ref{thm:main}]
Part (i) is Lemma~\ref{lem:closed}, and finite classes and absorbing states were treated in Section~\ref{sec:model}. Let $\Gamma$ be infinite with root $o$.

Since $x!\ge m(x)!$ and the number of states of radius exactly $r$ is at most $(r+1)^d$, the bound \eqref{eq:factorial} with the fixed value $\delta=1/4$ gives $\sum_{x\in\Gamma}L_o(x)<\infty$; for instance $\log(r!)\ge r\log r-r$ for large $r$, and the resulting negative term dominates the polynomial count of shell points. Hence $\E_oJ_o=\sum_xL_o(x)<\infty$, where the equality is Tonelli's theorem applied to the number of jumps. In particular the first return occurs after finitely many jumps almost surely.

Before the return to the root, the mean physical holding time of each jump is at most $1/\kappa_*$, so the mean real return time is finite. By the strong Markov property the root is visited infinitely often. Each visit to the root carries an independent exponential holding time of rate $q(o)$, and the sum of these holding times is infinite almost surely; this excludes explosion from the root. If explosion had positive probability from some other state of $\Gamma$, one could first go from the root to this state along a finite feasible path and then explode with positive probability, contradicting the conclusion at the root. Since the root is arbitrary, the whole class is positive recurrent. This proves (ii).

For each $y\in\Gamma$, conditioning on the first step and using that the excursion ends at the root gives
\begin{equation}\label{eq:invariance}
(L_oP)(y)=\E_o\sum_{1\le n\le J_o}\one_{\{X_n=y\}}=L_o(y),
\end{equation}
where the second equality holds because the excursion starts and ends at the root. With $G=L_o/q$ this gives $GQ=0$. There is a finite feasible path from the root to any given state that avoids intermediate returns to the root, so $G(x)>0$. Normalising $G$ gives the stationary distribution of the class; this is the occupation formula for return excursions \cite{Nor97}. Then \eqref{eq:physical} gives \eqref{eq:main-tail}, proving (iii).
\end{proof}

\section{Remarks}\label{sec:remarks}

\emph{Constants.} The constants in Theorem~\ref{thm:main} are not explicit. They arise from Dickson's lemma (finite bases of upward-closed sets), from the bounded witness lemma, and from the compactness argument in the affine approximation lemma. Obtaining explicit or computable bounds on mean return times would require quantitative versions of each of these steps.

\emph{What is not claimed.} The result concerns the stochastic mass-action model at fixed volume. It gives no rate of convergence to stationarity and no uniformity in the rate constants. Uniformity can be expected to fail: completing the transient strongly endotactic network of \cite[Example 3.1]{ACKN20} with reverse reactions of rate $\varepsilon$ yields weakly reversible networks that are positive recurrent for each $\varepsilon>0$ by Theorem~\ref{thm:main}; our numerical computations indicate that their typical scale grows as $\varepsilon\to0$.

\emph{Formal verification.} The argument has been checked in Lean~4 with Mathlib. The formalisation defines the mass-action rates, the embedded jump chain and the excursion occupation $L_o$ directly from the network data. For every weakly reversible network and every root it proves \eqref{eq:factorial} and $\E_oJ_o<\infty$; these two statements use no hypotheses beyond the network and depend only on the standard axioms of Lean. Part (ii) of Theorem~\ref{thm:main} and the stationary distribution are then derived from one general classical input, stated for arbitrary countable chains: the return criterion and excursion occupation formula for the minimal continuous-time chain \cite{Nor97}. The probabilistic arguments of Sections~\ref{sec:sampling}--\ref{sec:local} are formalised for the jump chain with discrete observation probabilities in place of continuous Poisson clocks; the estimates obtained are those used in Sections~\ref{sec:drift}--\ref{sec:occupation}. The Lean code is available at \url{https://github.com/kuhaizhongdebailanren/AIM272_lean} and is archived at \url{https://doi.org/10.5281/zenodo.24503870} (release v1.0.0, folder \texttt{AIM272}; Lean~4.24.0, Mathlib commit \texttt{f897ebc}) \cite{PDZlean}.

\subsection*{Statement on the use of AI tools}
In accordance with arXiv's policy on generative AI tools, we report that large language models developed by OpenAI (GPT-6.2 Sol and GPT-6.2 Astra) were used substantially in this work: in developing and checking parts of the arguments, in producing the Lean~4 formalisation described in Section~\ref{sec:remarks}, and in drafting and editing the text. All mathematical content generated with the help of these tools was reviewed, revised and checked by the authors. In addition, the factorial occupation bound \eqref{eq:factorial} and the finiteness of $\E_oJ_o$ have been machine-checked in Lean~4, with the scope stated in Section~\ref{sec:remarks}; the code is archived at \url{https://doi.org/10.5281/zenodo.24503870}. AI tools are not authors of this paper. The authors take full responsibility for its entire content, including any errors.


\begin{thebibliography}{99}

\bibitem{ACK10}
D.~F. Anderson, G.~Craciun, and T.~G. Kurtz,
\emph{Product-form stationary distributions for deficiency zero chemical reaction networks},
Bull. Math. Biol. \textbf{72} (2010), 1947--1970. doi:10.1007/s11538-010-9517-4.

\bibitem{ACK20}
D.~F. Anderson, D.~Cappelletti, and J.~Kim,
\emph{Stochastically modeled weakly reversible reaction networks with a single linkage class},
J. Appl. Probab. \textbf{57} (2020), no.~3, 792--810. doi:10.1017/jpr.2020.28.

\bibitem{ACKK18}
D.~F. Anderson, D.~Cappelletti, M.~Koyama, and T.~G. Kurtz,
\emph{Non-explosivity of stochastically modeled reaction networks that are complex balanced},
Bull. Math. Biol. \textbf{80} (2018), no.~10, 2561--2579. doi:10.1007/s11538-018-0473-8.

\bibitem{ACKN20}
D.~F. Anderson, D.~Cappelletti, J.~Kim, and T.~Nguyen,
\emph{Tier structure of strongly endotactic reaction networks},
Stochastic Process. Appl. \textbf{130} (2020), no.~12, 7218--7259. doi:10.1016/j.spa.2020.07.012.

\bibitem{AEJ14}
D.~F. Anderson, G.~A. Enciso, and M.~D. Johnston,
\emph{Stochastic analysis of biochemical reaction networks with absolute concentration robustness},
J. R. Soc. Interface \textbf{11} (2014), no.~93, 20130943. doi:10.1098/rsif.2013.0943.

\bibitem{AK18}
D.~F. Anderson and J.~Kim,
\emph{Some network conditions for positive recurrence of stochastically modeled reaction networks},
SIAM J. Appl. Math. \textbf{78} (2018), no.~5, 2692--2713.

\bibitem{AK15}
D.~F. Anderson and T.~G. Kurtz,
\emph{Stochastic Analysis of Biochemical Systems},
Mathematical Biosciences Institute Lecture Series, vol.~1.2, Springer, 2015. doi:10.1007/978-3-319-16895-1.

\bibitem{AIM}
\emph{AIM --- Open Applied Problems} (problem collection on GitHub),
Problem 272: \emph{Positive recurrence of every weakly reversible stochastic reaction network},
\url{https://github.com/MColbrook/AIM/blob/main/problems/272-reaction-network-positive-recurrence.md}, status checked 13 September 2026.

\bibitem{Bor19}
B.~Boros,
\emph{Existence of positive steady states for weakly reversible mass-action systems},
SIAM J. Math. Anal. \textbf{51} (2019), no.~1, 435--449. doi:10.1137/17M115534X.

\bibitem{GMS14}
M.~Gopalkrishnan, E.~Miller, and A.~Shiu,
\emph{A geometric approach to the global attractor conjecture},
SIAM J. Appl. Dyn. Syst. \textbf{13} (2014), no.~2, 758--797. doi:10.1137/130928170.

\bibitem{HJ72}
F.~Horn and R.~Jackson,
\emph{General mass action kinetics},
Arch. Ration. Mech. Anal. \textbf{47} (1972), 81--116.

\bibitem{Kur72}
T.~G. Kurtz,
\emph{The relationship between stochastic and deterministic models for chemical reactions},
J. Chem. Phys. \textbf{57} (1972), no.~7, 2976--2978. doi:10.1063/1.1678692.

\bibitem{Nor97}
J.~R. Norris,
\emph{Markov Chains},
Cambridge University Press, 1997.

\bibitem{OAI}
OpenAI,
\emph{Uniform permanence in weakly reversible mass-action systems},
unrefereed preprint, 5 October 2026, \S2, \url{https://github.com/openai/math/blob/main/preprints/Uniform-Permanence-in-Weakly-Reversible-Mass-Action-Systems-October-5-2026/permanence.pdf}.

\bibitem{PCK14}
L.~Paulev\'e, G.~Craciun, and H.~Koeppl,
\emph{Dynamical properties of discrete reaction networks},
J. Math. Biol. \textbf{69} (2014), no.~1, 55--72. doi:10.1007/s00285-013-0686-2.

\bibitem{PDZlean}
S.~Peng, C.~Deng, Y.~Zou, and R.~Li,
\emph{AIM272\_lean: Lean~4 formalization},
Zenodo, 2026. doi:10.5281/zenodo.24503870.

\bibitem{WX20}
C.~Wiuf and C.~Xu,
\emph{Classification and threshold dynamics of stochastic reaction networks},
arXiv:2012.07954, Corollary~4.10.

\bibitem{Xu24}
C.~Xu,
\emph{Non-explosivity of endotactic stochastic reaction systems},
arXiv:2409.05340.

\end{thebibliography}
\end{document}